\documentclass[review,hidelinks,onefignum,onetabnum]{siamart251104}
\usepackage[margin=0.75in]{geometry}

\input{ex_shared}

\ifpdf
\hypersetup{
  pdftitle={Representation of quasi-periodic functions and Hausdorff-Young inequalities for Besicovitch almost periodic functions},
  pdfauthor={Aihua Fan, Kai Jiang, and Pingwen Zhang}
}
\fi

\begin{document}

\maketitle
\begin{abstract}
 For a class of $\mathbb{R}^d$-, and $\mathbb{Z}^d$-actions on the $n$-dimensional torus $\mathbb{T}^n$, we characterize their unique ergodicity and establish a theorem of Weyl type.  This result allows us to establish an isomorphism between the Banach algebra of quasi-periodic functions with spectrum in a given $\mathbb{Z}$-module and the Banach algebra of periodic functions on a torus. This, in return, allows us to give a very simple proof of Hausdorff-Young inequalities for Besicovitch almost periodic functions. The regularity of the parent function of a quasi-periodic function is also studied. \\
\end{abstract}



\section{Introduction and main results}

We first study two dynamical systems acting on the $n$-dimensional torus $\mathbb{T}^n$ (with $n\ge 1$), one is a $\mathbb{R}^d$-action and the other is a $\mathbb{Z}^d$-action (with $d\ge 1$). 
Both are determined by a given set of $n$ vectors
$\bp_1, \cdots, \bp_n$ in the $d$-dimensional Euclidean space $\mathbb{R}^d$.
It is more practical for us to represent vectors in $\mathbb{R}^d$ as column vectors, i.e., $d\times 1$ matrix. Thus, the inner product $\bx\cdot \by$ is usually written as a matrix product $\bx^T \by$ (where $\bm{A}^T$ denotes the transpose of a matrix $\bm{A}$).  We  put the vectors $\bp_1, \cdots, \bp_n$  together as a $d\times n$ matrix
$$
\bP =(\bp_1, \bp_2,\cdots, \bp_n).
$$
Now, for any $\bx\in \mathbb{R}^d$, we define the map $\Phi_{\bx}: \mathbb{T}^n \to \mathbb{T}^n$ by 
$$
\Phi_{\bx}(\by) = (y_1 + \bp_1^T \bx , y_2+\bp_2^T \bx, \cdots, y_n + \bp_n^T \bx)^T ~\mbox{mod}~ \bbZ^n = (\by + \bP^T \bx)~\mbox{mod}~ \bbZ^n,
$$
where the $``+"$ in $y_j+\bp_j^T \bx$ is the operation in the group $\mathbb{T}$ and the $``+"$ in
$\by + \bP^T \bx$ is the operation in the group $\mathbb{T}^n$. 
It is clear that 
$$
\forall \bx, \bx'\in \mathbb{R}^d, \quad \Phi_{\bx}\circ  \Phi_{\bx'}=\Phi_{\bx+\bx'}
$$
Thus we get a topological dynamical system $(\mathbb{T}^n, (\Phi_{\bx})_{\bx\in \mathbb{R}^d})$. If we consider the restriction of $\bx$ on $\mathbb{Z}^d$, we get another dynamical system
$(\mathbb{T}^n, (\Phi_{\bx})_{\bx\in \mathbb{Z}^d})$.  
For these two dynamical systems, the normalized Lebesgue measure (Haar measure) is invariant. If the Lebesgue measure is the unique invariant measure, we say that the system is uniquely ergodic. We refer to \cite{walters1982introduction} for basic notions in dynamical systems and ergodic theory. 
The linear forms $\bx \mapsto \bp_j^T \bx$ are studied 
from the point of view of simultaneous Diophantine approximation. See Chapter III of \cite{Cassels1957}, where the Kronecker theorem is proved (p.53). 

We shall study the linear forms from a dynamical point of view and prove the following criterion for the unique ergodicity. It is a result stronger than the conclusion of the Kronecker theorem. We qualify this result as  
Kronecker-Weyl theorem. 

\begin{theorem}[Kronecker-Weyl for $\mathbb{R}^d$-action]\label{thm:UE1}
Consider the dynamical system
$(\mathbb{T}^n, (\Phi_{\bx})_{\bx\in \mathbb{R}^d})$. The following are equivalent.
\begin{enumerate}
\item[{\rm (a)}] The vectors $\bp_1, \cdots, \bp_n$ are $\bbQ$-independent, namely
\begin{equation}\label{eq:K-condQ}
\sum_{j=1}^n r_j \bp_j=0,\ \  r_j\in\bbQ \ (\forall j) \Longrightarrow r_1=r_2=\cdots=r_n=0.
\end{equation}
    \item[{\rm (b)}] The system is uniquely ergodic. 
    \item[{\rm (c)}] The Lebesgue measure is ergodic. 
    \item[{\rm (d)}]
    For any $F\in C(\mathbb{T}^n)$, we have
    \begin{equation}\label{eq:K-conl}
   \forall \by \in \mathbb{T}^n, \quad  \lim_{T\to \infty}\frac{1}{(2T)^d}\int_{[-T, T]^d} F(\Phi_{\bx}(\by)) \,d\bx =\int_{\mathbb{T}^n} F(\bz) \,d\bz
    \end{equation}
\end{enumerate}
\end{theorem}

Let us make the following remarks:
\begin{itemize}
    \item[{\rm (R1)}] 
That the $\mathbb{Q}$-independence in Theorem \ref{thm:UE1} (a) implies the ergodicity in Theorem  \ref{thm:UE1} (c) was  essentielly proved in \cite{JLZ2024}. Using Pitt's ergodic theorem (cf. \cite{Pitt1942}), the authors of \cite{JLZ2024} concluded that the equality in \eqref{eq:K-conl} holds true for almost all points $\by$ and for Lebesgue integrable functions $F$. We point out that this ensures the pointwise convergence as announced in Theorem \ref{thm:UE1} (d) for continuous functions $F$, if we remark that $\Phi_{\bx}$ is an isometry. Such results are fundamental for the projection method in scientific computation of quasi-periodic systems (cf. \cite{JLZ2024, JZ2014, JZ2018}). 

\item[{\rm (R2)}] It can be proved that the limit in \eqref{eq:K-conl} is uniform in $\by$. Because the $\mathbb{Q}$-independence implies the following decomposition
$$
C(\mathbb{T}^n)= \mathbb{C}\oplus \overline{{\rm Span} \left(\bigcup_{\bx\in \mathbb{R}^d} {\rm Im}\,(\Phi_{\bx} -\bI)\right)}.
$$
It is actually a common property for all unique ergodic amenable group actions.
See \cite{Fan2025} for details.  
\item[{\rm (R3)}] Actually, more general results have been proved in \cite{Fan2025} for a class of dynamical systems driven by locally compact Abelian groups ($\mathbb{R}^d$ and $\mathbb{T}^n$ are only two examples). 
\item[{\rm (R4)}] When the $\mathbb{Q}$-independence \eqref{eq:K-cond} is satisfied, the dynamics $(\mathbb{T}^n, \{\Phi_{\bx}\}_{\bx\in \mathbb{R}^d})$ is minimal, that is to say, for any $\by\in \mathbb{T}^n$, the orbit $\{\Phi_{\bx}(\by)\}_{\bx\in \mathbb{R}^d}$ is dense in $\mathbb{T}^n$.  It is a direct consequence of \eqref{eq:K-conl}: given any point $\by^*\in\mathbb{T}^n$
and fix any neighborhood of $\by^*$. Take $F\in C(\mathbb{T}^n)$
such that $0\le F(\by)\le 1$ for all $\by$,  with support in the neighborhood of $\by^*$ and $F(\by^*)=1$.  Then the limit in \eqref{eq:K-conl} shows that $\Phi_{\bx}(\by)$ falls into the neighborhood of $\by^*$ infinitely many times, because  the integral $\int_{\mathbb{T}^n} F(\by)\, d\by >0$.
\item[{\rm (R5)}] This minimality will be a very useful fact for our study on quasi-periodic functions, and the existence of the limit in (4) is another useful fact. See Lemma \ref{lem:lem1eq:norms} and Lemma \ref{lem:lem2} below.
   \end{itemize} 

For the $\mathbb{Z}^d$-action $(\mathbb{T}^n, (\Phi_x)_{x\in \mathbb{Z}^d})$, we have a similar result. This result will not be used for our study on almost periodic functions in the present paper.  But it is fundamental for the finite point recovery method in scientific computation of quasi-periodic systems \cite{jiang2024accurately}.

\begin{theorem}[Kronecker-Weyl for $\mathbb{Z}^d$-action]\label{thm:UE2}
Consider the dynamical system
$(\mathbb{T}^n, (\Phi_{\bx})_{\bx\in \mathbb{Z}^d})$. The following are equivalent.
\begin{enumerate}
\item[{\rm (a)}] The vectors $\bp_1, \cdots, \bp_n$ are $\mathbb{Z}$-independent $\mod \mathbb{Z}^d$, namely
\begin{equation}\label{eq:K-cond}
 \sum_{j=1}^n a_j \bp_j=0 \!\! \mod \mathbb{Z}^d, \ \ \ a_j \in \mathbb{Z} \ (\forall j) \Longrightarrow a_1=a_2=\cdots=a_n=0.
\end{equation}
    \item[{\rm (b)}] The system is uniquely ergodic. 
    \item[{\rm (c)}] The Lebesgue measure is ergodic. 
    \item[{\rm (d)}]
    For any $F\in C(\mathbb{T}^n)$, we have
    \begin{equation}\label{eq:K-conl2}
    \forall \by \in \mathbb{T}^n, \quad \lim_{T\to \infty}\frac{1}{(2T+1)^d}\sum_{\bx \in [-T, T]^d\cap \mathbb{Z}^d} F(\Phi_{\bx}(\by)) =\int_{\mathbb{T}^n} F(\bz)\, d\bz.
    \end{equation}
\end{enumerate}
\end{theorem}

We will apply Theorem \ref{thm:UE1} to study quasi-periodic functions. 
Let $AP(\mathbb{R}^d)$ be the space of Bohr
almost periodic functions (cf. \cite{LZ1982,Pankov1990}). 
Recall that $AP(\mathbb{R}^d)$ is a sub-Banach algebra of $C_b(\mathbb{R}^d)$ consisting of all uniform limits of trigonometric polynomials;   every $f$ in $AP(\mathbb{R}^d)$
 admits its Bohr mean
$$
\mathcal{M}(f):=\lim_{T\to \infty} \frac{1}{(2T)^d}
\int_{[-T, T]^d} f(\bx)\, d \bx 
$$
and its Bohr-Fourier spectrum $\sigma_f$ which is the countable set of  all
$\blam \in\mathbb{R}^n$ such that 
$\widehat{f}(\blam):=\mathcal{M}(f(\bx)e^{-2\pi i \blam^T \bx })\not=0$. We say that $f$ is quasiperiodic if the $\mathbb{Z}$-module generated  by $\sigma_f$ has an integral basis $\bp_1, \cdots, \bp_n$. For such a quasiperiodic function $f$, it can be proved that there exists a unique continuous periodic function $F\in C(\mathbb{T}^n)$ such that 
$f(\bx) =F(\bP^Tx)$ for all $\bx\in \mathbb{R}^d$ where $\bP^T$ is the transpose of the $d\times n$ matrix $\bP$ with column vectors $\bp_1, \cdots, \bp_n$. Let ${\rm QP}_{\bP}(\mathbb{R}^d)$ be the set of quasi-periodic functions
$f$ such that $\sigma_f\subset \bP(\mathbb{Z}^n)$. 
For $F \in L^1(\mathbb{T}^n)$ we denote by $\Sigma_F$ the Fourier spectrum of $F$ (the set of lattice points $\bk \in \mathbb{Z}^n$ such that $\widehat{F}(\bk)\not=0$). 

Let $W(\mathbb{T}^n)$ be the Wiener algebra of all continuous functions on $\mathbb{T}^n$ with absolutely summable Fourier coefficients (cf. \cite{Kahane1970}). This algebra is also denoted $A(\mathbb{T}^n)$.  We define 
$$
W^{QP}_{\bP}(\mathbb{R}^d)=\left\{f \in {\rm QP}_{\bP}(\mathbb{R^d}): 
\|f\|_{W^{QP}_{\bP}}:=\sum_{\blam \in \sigma_f} |\widehat{f}(\blam)|<\infty\right\}.
$$
We will see that $W^{QP}_P(\mathbb{R}^d)$ is a Banach algebra. We call it the Wiener algebra of quasi-periodic functions with spectrum in $\bP(\mathbb{Z}^n)$. 

Let $B^2(\mathbb{R}^d)$ be the space of Bescicovitch almost periodic functions of order $2$ on $\mathbb{R}^d$, which can be viewed as the restrictions on $\mathbb{R}^d$ of continuous functions on the Bohr compactification of $\mathbb{R}^d$ (see \cite{Pankov1990}, p.1-12). It is a Hilbert space equipped with the inner product $\mathcal{M}(f \overline {g})$
for $f, g \in B^2(\mathbb{R}^2)$ (such Bohr mean does exists) and $AP(\mathbb{R}^d)$ is continuously and densely embedded into $B^2(\mathbb{R}^d)$. Every $f\in B^2(\mathbb{R}^d)$ also admits a countable spectrum, also denoted $\sigma_f$. 
We denote by $B^2_{\bP}(\mathbb{R}^d)$ the space of all $f\in B^2(\mathbb{R}^d)$
having its spectrum $\sigma_f$ contained in $\bP\mathbb{Z}^n$. 

The following theorem shows that all three spaces ${\rm QP}_{\bP}(\mathbb{R}^d)$, $W^{QP}_{\bP}(\mathbb{R}^d)$ and $B^2_{\bP}(\mathbb{R}^d)$ are isomorphic to spaces of periodic functions.
In other words, such functions can be represented by periodic functions. 

\begin{theorem}\label{thm:QP}
   Suppose that  $\bP=(\bp_1, \bp_2, \cdots, \bp_n)$ is a $d\times n$ matrix such that ${\rm rank}_\mathbb{Q}\bP=n$. 
   Consider the mapping  $\mathcal{L}: C(\mathbb{T}^d) \to {\rm QP}_{\bP}(\mathbb{R}^d)$
   defined by $\mathcal{L}(F) = f$ where $f(\bx) = F(\bP^T \bx)$. Then
   \begin{itemize}
   \item[(a)]
   The mapping 
    is an isometric isomorphism of Banach algebras from $C(\mathbb{T}^d)$ onto ${\rm QP}_{\bP}(\mathbb{R}^d)$.
   \item[(b)] For $f(\bx) = F(\bP^T \bx)\in {\rm QP}_{\bP}(\mathbb{R}^d)$, we have $\sigma(f) = \bP(\Sigma_F)$, and for $\bk\in \Sigma_F$ we have
   $
          \widehat{f}(\bP\bk) = \widehat{F}(\bk).
   $
   \item[(c)] $W^{QP}_{\bP}(\mathbb{R}^d)$, equipped with the norm $\|f\|_{W^{QP}_{\bP}}$, is a Banach, which is isometrically isomorphic to the Wiener algebra $W(\mathbb{T}^n)$. The mapping $\mathcal{L}$ realizes the isomorphism.
   \item[(d)] If $f\in W^{QP}_{\bP}(\mathbb{R}^d)$ and $f(\bx)\not=0$ for all $\bx\in \mathbb{R}^d$, then $f^{-1} \in W^{QP}_{\bP}(\mathbb{R}^d)$.
   \item[(e)] $B^2_{\bP}(\mathbb{R}^2)$
   is a Hilbert space which  is isometrically isomorphic to $L^2(\mathbb{T}^n)$. 
   \end{itemize}
\end{theorem} 

\begin{figure}[!hbpt]
\centering
\includegraphics[width=11cm]{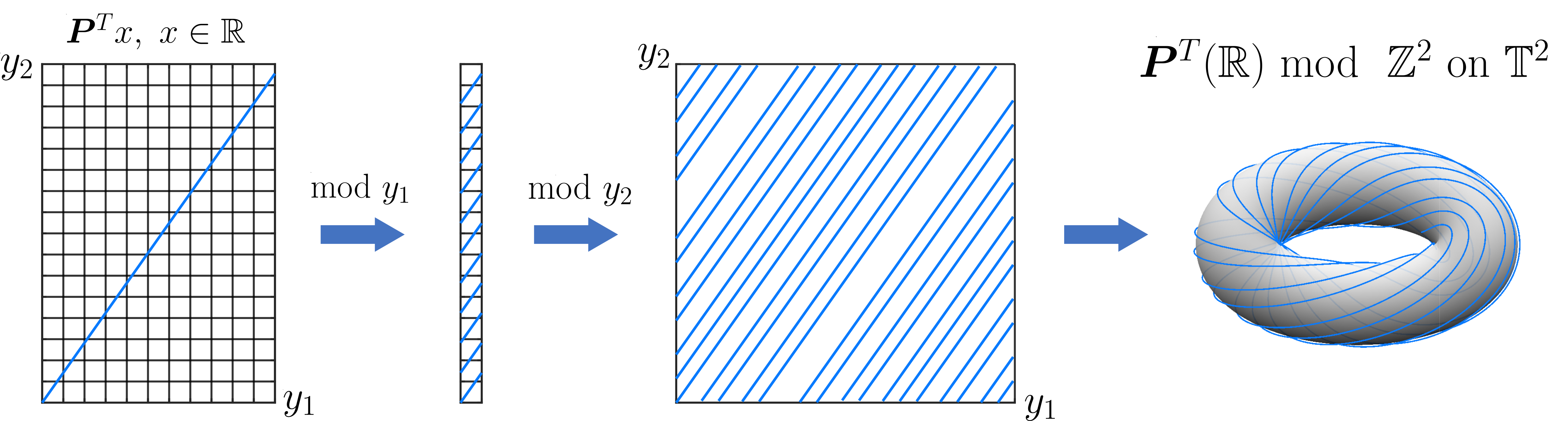}
\caption{One-dimensional irrational manifold $\bP^T (\bbR) \mod \bbZ^2 $ on $\bbT^2$ when $\bP = (1,\sqrt{2})$.}
\label{fig:pm_ergodic}
\end{figure}

When $f(\bx) =F(\bP^T \bx)$, the periodic function $F$ is called the parent function of $f$. From the point of view of theoretical study and scientific computation, the regularity of the parent is a key point. 
It is clear that if $F$ is of class $C^r$ for some $r\ge 1$ ($r$-times continuously differentiable), so is $f$. But the converse is far from clear and is not true in general.
The problem is that $f$ provides only information to $F$ on the $d$-dimensional submanifold $\bP^T(\mathbb{R}^n) \mod \mathbb{Z}^d$
of $\mathbb{T}^n$. In other words, there is no direct information about $F$ on $\mathbb{T}^n \setminus \bP^T(\mathbb{R}^d) \mod \mathbb{Z}^n$. 
Figure \ref{fig:pm_ergodic} presents an example of $\bP^T(\bbR)\mod \bbZ^2$ on $\bbT^2$, where $\bP=(1,\sqrt{2})$. 
Therefore there arises the question: how much regularity of $F$ can be observed from  that  of $f$?

This seems a tricky question. Let us give a first answer. 
Consider a quasi-periodic function defined on $\mathbb{R}$, which is of the form
\[
f(x) \sim  \sum_{\bk\in \mathbb{Z}^d}
a_{\bk} e^{2\pi i (\bP \bk) x}
\]
where $a_{\bk}$ is the Bohr-Fourier coefficient $\widehat{f}(\bP  \bk)$. Recall that the Bohr spectrum of $f$ is equal to 
\[
\sigma_f = \{ k_1 \bp_1 +\cdots + k_n \bp_n: a_{\bk} \neq 0 \},
\] 
which belongs to the subgroup  $p_1\mathbb{Z} + \dots + {p}_n\mathbb{Z}$ of $\mathbb{R}$ (it is also a $\mathbb{Z}$-module). It is possible that this subgroup is dense in $\mathbb{R}$. 
Indeed,  
 $p_1\mathbb{Z} + \dots + {p}_n\mathbb{Z}$ is dense in $\mathbb{R}$ if and only if $ p_i/p_j$ is irrational for some couple $i, j$ (see Lemma \ref{lem:density}). 
It is clear that 
$$
    |k_1{p}_1 +\cdots + k_n {p}_n| \le D \sqrt{k_1^2 +\cdots +k_n^2}
$$
for $D=\sqrt{p_1^2 +\cdots +p_n^2}$. We will make use of the following assumption that there exists $D'>0$ such that
\begin{equation}\label{eq:D-cond}
\forall \bk\in \sigma_f, \quad |k_1{p}_1 +\cdots + k_n {p}_n| \ge D' \sqrt{k_1^2 +\cdots +k_n^2}.
\end{equation}
This condition \eqref{eq:D-cond} is satisfied if $p_j>0$ for all $1\le j\le n$ and if $\sigma_f$ is contained in 
$$ \{p_1k_1+ \cdots + p_n k_n: (k_1, \cdots, k_n)\in \mathbb{N}^n\}.$$
So, the following theorem applies to quasi-periodic functions of the form
\[
f(x) = \sum_{k_1\ge 0,\cdots, k_n \geq 0} a_{k_1,\cdots, k_n} e^{2\pi i (p_1 k_1 + \cdots + p_n k_n) x}
\quad (p_1>0, \cdots, p_n >0).
\]

\begin{theorem}\label{thm:regularity}
Let \( f: \mathbb{R} \to \mathbb{C} \) be a quasi-periodic function with $p_1, \cdots, p_n$ as the integral basis of the module generated by its spectrum $\sigma_f$. Suppose that 
\begin{itemize}
    \item[(a)] there exist an integer $r>n$ and a real number $0<\eta<1$ such that  the $r$-times derivative $f^{(r)}$ is $\eta$-H\"{o}lder continuous, namely $f \in C^{r, \eta}(\mathbb{R})$;
\item[(b)] 
its Bohr spectrum $\sigma_f$ satisfies the discreteness condition \eqref{eq:D-cond}. 
\end{itemize}
Then the parent function $F: \mathbb{T}^n \to \mathbb{C}$ of $f$ is $(r-n)$-times continuously differentiable.
\end{theorem} 

We shall discuss Besicovitch-Sobolev spaces and prove a similar embedding theorem
for Besicovitch almost periodic functions of order $p$ with $1\le p<\infty$ (Theorem \ref{thm:SB-embedding}).

\medskip

The following Hausdorff-Young inequality was stated and proved by Avantaggiati et al for the case $d=1$ in \cite{ABI1995}. We are going to give a very simple proof of this fundamental result for the Besicovitch functions, which was a classical result in the case of periodic functions (cf. \cite{Katznelson2004}) 

\begin{theorem}[Hausdorff-Young inequality] \label{thm:HY} Let $1\le q<\infty$ with its conjugate
$q'=\frac{q}{q-1}$.
Assume $f \in B^q_{ap}(\mathbb{R}^d)$. 
\begin{itemize}
\item[(a)] If\ $1\le q\le 2$, we have
$$
    \left(
    \sum_{\blam \in \sigma_f} |\widehat{f}(\blam)|^{q'}
    \right)^{1/q'}
    \le \|f\|_q.
    $$
    \item[(b)] If \ $2\le q<\infty $, we have
$$
    \|f\|_q \le \left(
    \sum_{\blam\in \sigma_f} |\widehat{f}(\blam)|^{q'}
    \right)^{1/q'}.   
    $$
\end{itemize}
\end{theorem}


The key idea for proving Theorem \ref{thm:HY} is the observation that any trigonometric polynomial
is quasi-periodic. Then Theorem \ref{thm:QP} allows us to convert the problem to that for periodic functions.

\medskip The aforementioned theorems in will be proved sequence in the following sections.
\medskip

\section{Unique ergodicity of \texorpdfstring{$(\mathbb{T}^n, \{\Phi_{\mathbf{x}}\}_{{\mathbf{x}} \in \mathbb{R}^d})$}{}: Proof of Theorem \ref{thm:UE1}}
\label{sec:UE1}

Let $\mu$ be a probability Borel measure on $\mathbb{T}^n$. Notice that the Fourier coefficients of the image measure $\mu \circ \Phi_{\bx}^{-1}$ are easily computed:
$$
\forall \bk\in \mathbb{Z}^n, \qquad \widehat{\mu\circ \Phi_{\bx}^{-1}}(\bk)= \int_{\mathbb{T}^n} e^{-2\pi i \bk^T (\by+\bP^T \bx)}\,d \by = e^{-2\pi i \bk^T\bP^T \bx}\widehat{\mu}(\bk)
$$
This will be used.

(a) $\Longrightarrow$ (b). 
First we prove that the condition \eqref{eq:K-condQ} implies the unique ergodicity, i.e if $\mu$ is an invariant probability measure on $\mathbb{T}^n$, namely
$\mu \circ \Phi_{\bx}^{-1}=\mu$ for all $\bx\in \mathbb{R}^d$, then $\mu$ is the Lebesgue measure. Indeed, the invariance of $\mu$ can be described by Fourier coefficients of $\mu$ in the following manner
$$
\forall \bx\in \mathbb{R}^d, \qquad 
e^{-2\pi i (\bP \bk)^T \bx} \widehat{\mu}(\bk)=\widehat{\mu}(\bk).
$$
We claim that $\widehat{\mu}(\bk)=0$ for all $\bk\not= \bm{0}$. Otherwise, there exists $\bk^*\not=\bm{0}$ in $\mathbb{Z}^d$ such that $\widehat{\mu}(\bk^*)\not=0$ so that 
$$
\forall \bx\in \mathbb{R}^d, \qquad 
e^{-2\pi i (\bP \bk^*)^T \bx} = 1 \ {\rm or \ equivalently}\ 
(\bP \bk^*)^T \bx \in \mathbb{Z}.
$$
We can restate this fact as
$$
\mathbb{R}^d = \bigcup_{m\in \mathbb{Z}} H_m \quad {\rm with} \quad 
H_m :=\{\bx\in \mathbb{R}^d: (\bP \bk^*)^T \bx  = m\}.
$$
Since $\bk^*\not=\bm{0}$, the condition \eqref{eq:K-condQ} implies the non nullity of the vector $\bP\bk^*$:
$$
\bP \bk^* = k^*_1\bp_1
+ k^*_2\bp_2+\cdots  +k^*_n\bp_n\not= \bm{0}.
$$
Thus, each $H_m$ is a hyperplane and then has dimension $d-1$. Consequently, $\mathbb{R}^d$ would have dimension $d$, a contradiction. 

(b) $\Longrightarrow$ (c). It is trivially true.

(c) $\Longrightarrow$ (a). We prove it by contradiction. Suppose that \eqref{eq:K-condQ} is false. It suffices to exhibit a non-zero and non-constant Lebesgue integrable invariant function. Suppose, without loss of generality, that
$$
\bp_n =r_1\bp_1 + \cdots +r_{n-1}\bp_{n-1},\quad (r_j \in \mathbb{Q}).
$$
Let $b$ be the least common multiple of the denominators 
of $r_j$'s. We can write the above equality as
\begin{equation}\label{eq:non-ergodic}
b\bp_n =a_1\bp_1 + \cdots +a_{n-1}\bp_{n-1},\quad (a_j \in \mathbb{Z}).
\end{equation}
Consider the non-trivial group character of $\mathbb{T}^n$ defined by
$$
\gamma(\by)=\exp\Big[2\pi i \big(a_1y_1 +\cdots +a_{n-1}y_{n-1} - b y_n\big)\Big].
$$
We check that $\gamma(\cdot)$ is a desired invariant function. Indeed,
\begin{eqnarray*}
\forall \bx\in \mathbb{R}^d, \quad 
\gamma \circ \Phi_{\bx}(\by) 
&=& \gamma (\by+\bP^T \bx)
= \gamma(\by) \gamma(\bp_1^T \bx, \cdots, \bP_n^T \bx)=\gamma(\by)
\end{eqnarray*}
because
$$
\gamma(\bp_1^T \bx, \cdots, \bp_n^T \bx)=
\exp\Big[2\pi i \big(a_1\bp_1^T \bx +\cdots +a_{n-1}\bp_{n-1}^T \bx - b\bp_n^T \bx \big)\Big]=1,
$$
where the last equality results from \eqref{eq:non-ergodic}.

(d) $\Longrightarrow$ (b). For simplicity, we denote by $\langle \mu, F\rangle$ the integral of $F$ with respect to a measure $\mu$ and we denote by $\nu$ the Lebesgue measure on $\mathbb{T}^d$. Assume that $\mu$ is a probability invariant measure.  Let
$$
A_T F(\by) =\frac{1}{(2T)^d}\int_{[-T, T]^d} F\circ\Phi_{\bx}(\by))\, d\bx. 
$$
We are going to prove that $\mu =\nu$, by showing 
$\langle \mu, F\rangle = \langle \nu, F\rangle $ for all $F\in C(\mathbb{T}^n)$.
Indeed,  as $\langle \mu, F\circ \Phi_x\rangle = \langle \mu, F\rangle$ for all $\bx\in \mathbb{R}^d$, by Fubini theorem and Lebesgue's dominated convergence theorem we have
$$
\langle \mu, F\rangle = 
\frac{1}{(2T)^d}\int_{[-T, T]^d} \langle \mu, F\circ \Phi_{\bx}(\cdot))\, d\bx=\langle \mu, A_TF(\cdot)\rangle \to \langle \nu, F\rangle,
$$
as $T\to +\infty$, 
where we have used the pointwise convergence of $A_TF(y)$ for the above limit. 

(b) $\Longrightarrow$ (d). As (b) and (c) are equivalent, we can use the ergodicity of the Lebesgue measure to ensure that 
$A_TF(y)$ converges almost everywhere to $\langle \nu, F\rangle$. This implies the everywhere convergence because $\Phi_{\bx}$ are isometries. Indeed, take a dense set $D$ in $\mathbb{T}^n$ such that for any $\by\in D$, $A_TF(\by)$ tends to $\langle \nu, F\rangle$.  For any $\epsilon >0$,
there exists $\delta >0$ such that 
$$
\by', \by'' \in \mathbb{T}^n, \ \|\by'-\by''\|<\delta \Longrightarrow |F(\by')-F(\by'')|<\epsilon.
$$
For any $\by\in \mathbb{T}^n$, take a $\by^*\in D$ such that $\|\by-\by^*\|<\delta$. The isometry implies
$$
\forall \bx\in \mathbb{R}^d, \quad |F\circ\Phi_{\bx}(\by)-F\circ \Phi_{\bx}(\by^*)| <\epsilon.
$$
Therefore
\begin{eqnarray*}
|A_TF(\by)-\langle \nu, F\rangle|&\le& |A_TF(\by)- A_TF(\by^*)| +|A_TF(\by^*)-\langle \nu, F\rangle| \\
&\le & \epsilon + |A_TF(\by^*)-\langle \nu, F\rangle|.
\end{eqnarray*}
The last term can be made smaller that $2\epsilon$ when $T$ is large. 
We have thus proved that $A_TF(\by)$ tends to $ \langle \nu, F\rangle$ for all $\by$, and uniformly in $\by$.

\section{Unique ergodicity of \texorpdfstring{$(\mathbb{T}^n, (\Phi_{\mathbf{x}})_{\mathbf{x} \in \mathbb{Z}^d})$}{}: Proof of Theorem \ref{thm:UE2}}
\label{sec:UE2}

The idea of proof is the same as that of last theorem. Let us sketch the proof by pointing out the different points.

    (a) $\Longrightarrow$ (b). 
Suppose that $\mu$  an invariant probability measure, namely
$$
\forall \bx\in \mathbb{Z}^d, \qquad 
e^{-2\pi i (\bP \bk)^T \bx} \widehat{\mu}(\bk)=\widehat{\mu}(\bk).
$$
Here $x$ runs over $\mathbb{Z}^d$, not $\mathbb{R}^d$. 
It suffices to prove  that $\widehat{\mu}(\bk)=0$ for all $\bk\not=\bm{0}$. It it were not the case, there would exist $\bk^*\not=\bm{0}$ in $\mathbb{Z}^d$ such that $\widehat{\mu}(\bk^*)\not=0$ so that 
$$
\forall \bx\in \mathbb{Z}^d, \qquad 
e^{-2\pi i (\bP \bk^*)^T \bx} = 1 \ {\rm or \ equivalently}\ 
(\bP \bk^*)^T \bx \in \mathbb{Z}.
$$
This is equivalent to 
$$
(\bP \bk^*)^T \be_j \in \mathbb{Z} \ \ {\rm for\ all}\ \ 1\le j \le n
$$
where $\{\be_1, \cdots, \be_n\}$ is the canonical basis of the lattice $\mathbb{Z}^d$.
In other words, each coordinate of the vector $\bP \bk^*=\sum_{j=1}^n k_j^* \bp_j$ is an integer, namely
$$\sum_{j=1}^n k_j^* \bp_j= \bm{0} \mod \mathbb{Z}^d.$$ 
Then, using the condition (a), we get $\bk^*=\bm{0}$, a contradiction.

(b) $\Longrightarrow$ (c). It is trivially true.

(c) $\Longrightarrow$ (a). We prove it by contradiction. Suppose that (a) is false. 
Then there exists a non-zero  integral vector $(a_1^*, \cdots, a_n^*)\in \mathbb{Z}^n\setminus\{\bm{0}\}$ such that $  \sum_{j=1}^n a_j^* \bp_j \in \mathbb{Z}^d$. 
We can check that the non-trivial group character defined by
$$
\gamma(\by)=\exp\Big[2\pi i \big(a_1^*y_1 +a_{2}^*y_{2}+\cdots +a_n^* y_n\big)\Big]
$$
is a (non-constant) invariant function, which will contradicts the ergodicity of the Lebesgue measure.  Indeed,
\begin{eqnarray*}
\forall \bx\in \mathbb{Z}^d, \quad 
\gamma \circ \Phi_{\bx}(\by) 
&=& \gamma (\by+\bP^T \bx)
= \gamma(\by) \gamma(\bp_1^T \bx, \cdots, \bp_n^T \bx)=\gamma(\by)
\end{eqnarray*}
where the last equality is due to
$$
\gamma(\bp_1^T \bx, \cdots, \bp_n^T \bx)=
\exp\Big[2\pi i \big(a_1^*\bp_1^T  +a_{2}^*\bp_{2}^T + \cdots + a_n^*\bp_n\big) \bx\Big]=1.
$$
This is because both vectors $a_1^*\bp_1^T  +a_{2}^*\bp_{2}^T + \cdots + a_n^*\bp_n$ (row) and $\bx$ (column) are integral   and then their (inner) product is an integer. 

The proof of other points are the same.

\section{Representations of quasi-periodic functions: Proof of Theorem \ref{thm:QP}}
\label{sect:QP}

The proof of Theorem \ref{thm:QP} is based on  the following two lemmas. Both of them are corollaries of Theorem \ref{thm:UE1}. 

\begin{lemma}\label{lem:lem1eq:norms}
    Let $f(\bx)=F(\bP^T \bx)$
    for some $F\in C(\mathbb{T}^n)$. Suppose that $\bP$ is $d\times n$ matrix with ${\rm rank }_{\bbQ} \bP=n$. Then 
    \begin{equation}\label{eq:norms}
     \|f\|_\infty =\|F\|_\infty,    
    \end{equation}
    where
    $
    \|f\|_\infty =\sup_{\bx\in \mathbb{R}^d} |f(\bx)|, 
    ~ \|F\|_\infty =\sup_{\by\in \mathbb{T}^n} F(\by)|.
    $
\end{lemma}
Indeed, this unique ergodicity implies the minimality of the flow.
The equality $\|f\|_\infty=\|F\|_\infty$ is a consequence of   the fact that $\{\bP^T \bx\}_{\bx\in \mathbb{R}^d}$, the orbit
of the neutral point of $\mathbb{T}^n$,  is dense in 
$\mathbb{T}^n$ according to the the minimality. 

\begin{lemma}\label{lem:lem2}
    Under the same condition as in Lemma \ref{lem:lem1eq:norms},  we have
\begin{equation}\label{eq:QP}
\mathcal{M}(f):=\lim_{T\to \infty} \frac{1}{(2T)^d}
\int_{[-T, T]^d} f(\bx)\, d\bx =\int_{\mathbb{T}^n} F(\bz)\, d\bz.
\end{equation}
\end{lemma}
Indeed, the limit in \eqref{eq:QP} is nothing but the time average $F$ along the orbit of the neutral element of $\mathbb{T}^n$ which is equal to the space mean value of $F$, by the unique ergodicity.  
\medskip

The proofs of Theorem \ref{thm:QP} (a) and (b) consist of the following four steps. We will use the above two lemmes, Bochner's criterion for almost periodicity and Parseval's identity for almost periodic functions.  
 \medskip 

 Step 1. {\it Suppose $F\in C(\mathbb{T}^n)$. Then $f(\bx)$ defined by $F(\bP^T \bx)$ is almost periodic, $\bP(\Sigma_F) \subset \sigma_f$ and }
 \begin{equation}\label{eq:FF}
     \forall \bk \in \mathbb{Z}^n, \ \ \   \widehat{f}(\bP \bk) = \widehat{F}(\bk).
\end{equation}

As $F$ is periodic, the set of translates $\{F(\cdot +\by)\}_{\by\in \mathbb{T}^n}$
is relatively compact in $C(\mathbb{T}^n)$ by the Arzela-Ascoli theorem. Then according to Bochner's criterion, the set of
translates $\{f(\cdot + \bx)\}_{\bx\in \mathbb{R}^d}$ of $f$ is relatively compact  under supremum norm, by Lemma \eqref{lem:lem1eq:norms}. So $f$ is almost periodic, by the Bochner criterion. Now, by Lemma \eqref{lem:lem2}, we compute
$$
  \widehat{f}(\bP\bk)
  = \mathcal{M}(F(\bP^T\bx)e^{-2\pi i (\bP \bk)^T \bx})
  =\mathcal{M}(F(\bP^T\bx)e^{-2\pi i  \bk^T \bP^T \bx}) =\widehat{F}(\bk).
$$
We have thus proved \eqref{eq:FF}. This equality implies  $\bP(\Sigma_F)\subset \sigma_f$. But for the moment we don't confirm that $\bP(\Sigma_F)\supset \sigma_f$, because the above computation is performed only for $\bP\bk$'s. 
\medskip 

Step 2. {\it We have $\sigma_f=\bP(\Sigma_F)$ for 
$f(\bx)$ defined by $F(\bP^T \bx)$, so $f$ belongs to ${\rm QP}_{\bP}(\mathbb{R}^d)$.
}\\
By Parseval's relation for almost periodic functions (cf.~\cite{LZ1982}, p.\,31) and Parseval's relation for periodic functions, we have
$$
      \mathcal{M}(|f|^2) = \sum_{\blam \in \sigma_f} |\widehat{f}(\blam)|^2, \qquad \|F\|_{L^2}^2 = \sum_{\bk \in \Sigma_F} |\widehat{F}(\bk)|^2.
$$ 
By \eqref{eq:QP} applied to $|F|^2$, we have $\mathcal{M}(|f|^2) =  \|F\|_{L^2}^2$. This is a key point. Therefore 
$$
 \sum_{\blam \in \sigma_f} |\widehat{f}(\blam)|^2= \sum_{\bk\in \Sigma_F} |\widehat{F}(\bk)|^2 =  \sum_{\blam \in \bP(\Sigma_F)} |\widehat{f}(\blam)|^2
 $$
 where \eqref{eq:FF} is used once more for the last equality. As  $\bP(\Sigma_F)\subset \sigma_f$,
 this implies  $\sigma_f \setminus \bP(\Sigma_F)$ must be empty. Thus we have proved $\sigma_f =\bP(\Sigma_F)$.

 Step 3. {\it The mapping  $\mathcal{L}: C(\mathbb{T}^n) \to {\rm QP}_{\bP}(\mathbb{R}^d)$
   defined by $F \mapsto f(\bx) = F(\bP^T \bx)$ is bijective.
   }

   The injectivity for the linear mapping $\mathcal{L}$
   is an immediate consequence of \eqref{eq:norms}. 
   We have only to prove the surjectivity. 
   Given a function $f\in {\rm QP}_{\bP}(\mathbb{R}^d)$ with its Bohr-Fourier series
   $$
     f(\bx) \sim \sum_{\bk\in \mathbb{Z}^n} \widehat{f} (\bP \bk) e^{2\pi i (\bP \bk)^T \bx }.
   $$
   Our aim is to find a continuous function 
   $F\in C(\mathbb{T}^n)$ such that $F(\bP^T \bx)=f(\bx)$. The data of $f$ only provides information about $F$ on the dense set $\bP(\mathbb{R}^d) \mod \mathbb{Z}^n$ of $\mathbb{T}^n$. But this information is sufficient to determine $F$. Our candidate is the function defined by the formal series
   $$
      F(\by) \sim \sum_{\bk\in \mathbb{Z}^n} \widehat{f} (\bP \bk) e^{2\pi i \bk^T \by }.
   $$
   We are going to show that this function $F$ is continuous by showing that the F\'ejer sums of  this series uniformly converge. Recall that the
   $N$-th order F\'ejer sum of $F$ is defined by
   $$
     \sigma_N(F)(\by) =\sum_{|\bk|\le N}
     \left(1- \frac{|k_1|}{N}\right)
     \cdots \left(1- \frac{|k_n|}{N}\right) \widehat{f} (\bP \bk)  e^{2\pi i \bk^T \by}
   $$
   where $\bk =(k_1, \cdots, k_n)\in \mathbb{Z}^n$, $|\bk| =\max_{1\le j\le n} |k_j|$. Observe that
   \begin{equation}\label{eq:SS}
     \sigma_N(F)(\bP^T \bx) =\sum_{|\bk|\le N}
     \left(1- \frac{|k_1|}{N}\right)
     \cdots \left(1- \frac{|k_n|}{N}\right) \widehat{f} (\bP \bk)  e^{2\pi i (\bP\bk)^T \bx}=: \tilde{\sigma}_N(f)(\bx)
   \end{equation}
   where $\tilde{\sigma}_N(f)(\bx)$ is a Bochner-F\'ejer sum of $f$.
   This equality \eqref{eq:SS}, together \eqref{eq:norms}, implies
   $$
   \| \sigma_N(F) - \sigma_{N+p}(F) \|_\infty
   = \|\tilde{\sigma}_N(f) - \tilde{\sigma}_{N+p}(f) \|_\infty
   $$
   As Bochner-F\'ejer theorem asserts (cf. \cite{LZ1982}, p. ),  $\tilde{\sigma}_N(f)$ converges uniformly to $f$ and is thus a Cauchy sequence. Therefore $\sigma_N(F)$
   is a Cauchy sequence in $C(\mathbb{T}^n)$. Then
   $F$ is continuous and 
   we get $F(\bP^T \bx)=f(\bx)$
   for $\bx\in \mathbb{R}^d$ from \eqref{eq:SS}.
   \medskip

Step 4. {\it $\mathcal{L}: C(\mathbb{T}^d) \to {\rm QP}_{\bP}(\mathbb{R}^d)$ 
    is an isometric isomorphism of algebras}.\\
    We have seen that $\mathcal{L}$ is bijective. It is clearly an isomorphism of algebra. It is also isometric because of Lemma \ref{lem:lem1eq:norms}.
    \medskip
   
(c) is an immediate consequence of (b).

(d) is a consequence of (c) and the famous Wiener theorem for the algebra $A(\mathbb{T}^n)$. 

(e) Define $\mathcal{L}: L^2(\mathbb{T}^n) \to B^2_{\bP}(\mathbb{R}^d)$ by
\begin{equation}\label{eq:L2}
\mathcal{L}(F)(x) =\sum_{\bk \in \mathbb{Z}^n}
\widehat{F}(\bk) e^{2\pi i (\bP \bk)^T  \bx}.
\end{equation}
It is easy to check that $\mathcal{L}$ is the desired isomorphism. Thus we have proved Theorem \ref{thm:QP}.
\medskip

Remark that the mapping $F \mapsto f(\bx) :=F(\bP^T \bx)$ from $C(\mathbb{T}^n)$ onto ${\rm QP}_{\bP}(\mathbb{R}^d)$ is the restriction of the  mapping $\mathcal{L}$
 defined by \eqref{eq:L2}.

\section{Hausdorff-Young inequality in Besicovitch spaces \texorpdfstring{$B^p_{ap}(\mathbb{R}^d)$}{}}
\label{sec:HYineqn}

In this section, we are going to give a very simple proof of the Hausdorff-Young inequality for Besicovitch almost periodic functions, which was proved by Avantaggiati et al \cite{ABI1995} when $d=1$. Our proof is based on the observation that every trigonometric polynomial is quasi-periodic and it can be represented by its parent function (cf.~Theorem \ref{thm:QP}).   The the Hausdorff-Young inequality is a fundamental result. 
We will use it to study the regularity of parent functions in the next section. 

\subsection{Space \texorpdfstring{$B^q_{ap}(\mathbb{R}^d)$}{} of almost periodic functions in the sense of Besicovitch}

Let $1\le q<\infty$. The space $\mathcal{T}(\mathbb{R}^d)$ of trigonometric polynomials on $\mathbb{R}^d$ is equipped with the following norm
$$
\|S\|_q := \left(\mathcal{M}(|S|^q) \right)^{1/q}=:\lim_{T\to\infty}\left(\frac{1}{(2T)^d}\int_{[-T, T]^d} |S(\bx)|^q \, d\bx\right)^{1/q}.
$$
The above limit does exist, because $|S|^q$ is almost periodic in  Bohr sense.
The space $B^q_{ap}(\mathbb{R}^d)$ (or simply denoted $B^q_{ap}$) of Besicovitch almost periodic functions of order $q$ is defined to be the completion of $\mathcal{T}(\mathbb{R}^d)$ with respect to the norm $\|\cdot\|_q$ defined above (cf.~\cite{Besicovitch1932}). Attention: $\|\cdot\|_q$ is not the usual norm of Lebesgue space $L^q(\mathbb{R}^d)$.  Formally $B^\infty_{ap}(\mathbb{R}^d)$ is similarly defined with the norm $\|\cdot\|_\infty$ (the usual norm of supremum) and  we get 
$B^\infty_{ap}(\mathbb{R}^d)=AP(\mathbb{R}^d)$, the space of almost periodic functions in Bohr sense. 
We have the following obvious continuous embedding
$$
B^\infty_{ap} \subset B^{q_2}_{ap} \subset B^{q_1}_{ap}\subset  B^1_{ap} \quad (1<q_1<q_2<\infty)
$$
and $B^1_{ap}(\mathbb{R}^d)$
is the largest among all these spaces. 

For any polynomial $S$, the Fourier coefficient $\widehat{S}(\blam)=\mathcal{M}(S(\bx) e^{-i\blam^T \bx})$ is well defined for every $\blam \in \mathbb{R}^d$. The following obvious uniform estimate will be useful
\begin{equation}\label{eq:Coef-Estimate}
|\widehat{S}(\blam)-\widehat{Q}(\blam)|\le \|S-Q\|_1 \quad (\forall S, Q\in \mathcal{T}(\mathbb{R}^d)).
\end{equation}
By definition, an element $f$ of $B^1_{ap}(\mathbb{R}^d)$ is the limit of Cauchy 
sequence $f=(S_n)$ of polynomials. The inequality
\eqref{eq:Coef-Estimate} implies that $(\widehat{S}_n(\blam)$ is a Cauchy sequence, so that the following limit exists
\begin{equation}\label{eq:FourierCoef-B}
\forall \blam \in \mathbb{R}^d, \quad 
\widehat{f}(\blam) :=\lim_{n\to \infty} \widehat{S}_n(\blam).
\end{equation}
This defines the Fourier coefficient of $f$ at $\blam$. The spectrum of $f\in B^1_{ap}(\mathbb{R}^d)$ is defined as
$$
\sigma_f=\{\blam\in \mathbb{R}^d: \widehat{f}(\blam)\not=0\}.
$$
The spectrum $\sigma_f$ is at most countable for every $f\in B^1_{ap}(\mathbb{R}^d)$ (cf.~\cite{Katznelson2004}). Indeed
$$
\sigma_f \subset \varlimsup_{n\to \infty} \sigma_{S_n}.
$$
Thus each $f\in B^1_{ap}(\mathbb{R}^d)$ is associated to its Fourier series (or Fourier-Besicovitch series)
$$
f(\bx) \sim \sum_{\blam \in \sigma_f} \widehat{f}(\blam) e^{i\blam^T \bx}.
$$

For any Besicovitch function $f\in B^1_{ap}(\mathbb{R}^d$ and any polynomial $S$, the following quantity is well defined
$$
(f, S):= \mathcal{M}(f\overline{S}) = \sum \widehat{f}(\blam) \overline{\widehat{S}(\blam)}.
$$
It can be considered as an inner product if $f\in B^2_{ap}(\mathbb{R}^d)$. In general case, we can consider it as a duality. Actually the norm $\|f\|_q$ for $f\in B^q_{ap}(\mathbb{R}^d)$ can be represented by this duality through polynomials
\begin{equation}\label{eq:normdual}
    \|f\|_q = \sup_{S\in \mathcal{T}(\mathbb{R}^d), \|S\|_{q'}\le 1} |(f, S)|.
\end{equation}

\subsection{Proof of Hausdorff-Young theorem on \texorpdfstring{$B^q_{ap}(\mathbb{R}^d)$}{}}

Let us first prove these inequalities for polynomials using the classical Hausdorff-Young inequalities. This is possible just because every polynomial is quasi-periodic and every quasi-periodic function is well determined by its periodic parent function to which we can apply the classical Hausdorff-Young inequality.

Assume that  $f$ is a polynomial. Let $F \in C(\mathbb{T}^n)$ be its parent function and $f(\bx)=F(\bP^T \bx)$  with projection matrix $\bP$.
Let $\Sigma_F$ be the spectrum of $F$. We have 
$\sigma_f =\bP(\Sigma_F)$
and $\widehat{F}(\bk) =\widehat{f}(\bP \bk)$ (cf. Theorem \ref{thm:QP}). 

Assume $1\le q\le 2$. By the classical Hausdorff-Young inequality applied to $F$, we get
$$
\|f\|_q|= \|F\|_q \ge 
\left(
    \sum_{\bk\in \Sigma_F} |\widehat{F}(\bk)|^{q'}
    \right)^{1/q'}
    =  \left(
    \sum_{\bk\in \Sigma_F} |\widehat{f}(\bP \bk)|^{q'}
    \right)^{1/q'}
    =\left(
    \sum_{\blam\in \sigma_f} |\widehat{f}(\blam)|^{q'}
    \right)^{1/q'},
$$
where the first equality is a consequence of the identity\eqref{eq:QP} applied to $|f|^q$ in the place of $f$. This is the key point of the present proof. We have thus proved the first inequality for polynomilas. The second inequality for the case $q\le 2<\infty$ can be proved in the same way.

Now let $f$ be an arbitrary function in $B^q_{ap}(\mathbb{R}^d)$, with $1\le q\le 2$. Let $(f_n)$ be a Cauchy sequence tending to $f$. By what we have just proved we have
$$
\left(
    \sum_{\blam\in \sigma_{f_n}} |\widehat{f}_n(\blam)|^{q'}
    \right)^{1/q'}
    \le \|P_n\|_q.
$$
Fix an arbitary finite subset $F\subset \sigma_f$. 
Taking limit allows us to get
$$
\left(
    \sum_{\blam\in F }|\widehat{f}(\blam)|^{q'}
    \right)^{1/q'}
    \le \|f\|_q.
$$
Since $F$ is arbitrary, we have thus finished the proof of the first inequality.

To prove the second inequality, we will use the first one and a dual argument. For any trigonometric polynomial $Q$, we have
$$
(f_n, Q) = \sum \widehat{f_n}(\blam)  \overline{\widehat{Q}(\blam)}
$$
This is s finite sum, a sum over the spectrum of $Q$. It follows that
$$
|(f, Q)|=\lim_{n\to \infty}|
(f_n, Q)| = \left|\sum \widehat{f}(\blam)  \overline{\widehat{Q}(\blam)}\right|\le \left(\sum |\widehat{f}(\blam)|^{q'} \right)^{1/q'} \left(\sum |\widehat{Q}(\blam)|^{q} \right)^{1/q}
$$
where we have used \eqref{eq:FourierCoef-B} and the H\"{o}lder ineqaulity.
Then, 
by the Hausdorff-Young inequality already proved, we obtain
$$
|(f, Q)|\le 
\left(\sum_{\blam\in \Sigma_f} |\widehat{f}(\blam)|^{q'} \right)^{1/q'} \|Q\|_{q'}.
$$
We conclude by using the formula \eqref{eq:normdual}.

\section{Regularity of  parent functions}\label{sect:regularity}

Recall the relation between a quasi-periodic function $f$ on $\mathbb{R}^d$ and its parent
function $F$ on $\mathbb{T}^n$:
\begin{equation}\label{eq:f=F}
f(\bx) = F(\bP^T x)
\end{equation}
where $\bx \in \mathbb{R}^d$, $F \in C(\mathbb{T}^n)$, $\bP = [\bp_1, \ldots, \bp_n]$ is a $d \times n$ matrix with linearly $\mathbb{Q}$-independent columns $\bp_1, \ldots, \bp_n$. Actually the set $\{\bp_1, \ldots, \bp_n\}$ is an integral basis of the $\mathbb{Z}$-module generated by the Bohr spectrum of $f$.

If the parent function $F$ is continuously differentiable, so is the quasi-periodic function $f$. But the converse is far from true. How much regularity of
the parent function $F$ can be deduced from the regularity of the parent function $f$ ?

\subsection{Preliminary discussion}

Let $\by = \bP^T \bx$ be the point on the submanifold $\bP^T(\mathbb{R}^d)$ associate to $\bx$. 
From \eqref{eq:f=F} we get immediately
\[
 D_{\bP^T \bu} F(\by)= D_{\bu} f(\bx).
\]
That means, the directional derivative 
$D_{\bP^T\!\bu} F(\by)$ in the direction 
$\bP^T\!\bu$ of $F$ exists if and only if the directional derivative $D_{\bu} f(\bx)$ exits, and in this case they are equal.
Therefore, if we assume that $f$ is differentiable at $\bx$, then the directional derivatives $D_{\bP^T\! \be_i} F(\by)$  $(1 \leq i \leq d)$ exist where  $\{\be_1, \dots, \be_d\}$ is the canonical basis of $\mathbb{R}^d$, and
\[
  D_{\bP^T \be_i} F(\by) =\frac{\partial f}{\partial x_i} (\bx) \quad (1 \leq i \leq d).
\]
The directions $\bP^T\! \be_1.\cdots, \bP^T \! \be_d$ are independent. But there is a lack of $n-d$ directional derivatives, if we would like to  deduce the differentiability of $F$ at $\by$.

We are going to see that there exists quasi-periodic function of class $C^\infty$, but its parent function is not differentiable. One of reasons for such pathology is that the spectrum of the quasi periodic function has accumulation points. 

We will study the regularity of the parent function, basing on that of the corresponding quasi-periodic function, by using 
 the Fourier method. Basic ideas are based
on the following facts:
\begin{enumerate}
    \item $\widehat{F}(\bk) = \widehat{f}(\bP \bk)$, for all $\bk \in \mathbb{Z}^n$ (cf. Theorem \ref{thm:QP})
    \item Quick decay of $\widehat{F}$ implies the regularity of $F$ (it is classical, cf. \cite{Katznelson2004})
    \item The regularity of $f$ implies the decay of $\widehat{f}$, which in turn implies the decay of $\widehat{F}$.
\end{enumerate}
The fact (3) is precised in the following lemma. 
Let us first introduce the {\it Bohr $L^1$-modulus of continuity}  of an almost periodic function $f$ by
$$
  \omega_1(f, \delta): = \sup_{|\bt|\le \delta}\mathcal{M}(|f(\cdot + \bt)- f(\cdot)|).
$$
It is clear that $\omega_1(f, \delta) = O(\delta^\eta)$
if $f$ is $\eta$-H\"older in the sense that there exists a constant $C>0$ such that 
\[
\forall \bx, \bt \in \mathbb{R}^d, \quad |f(\bx+\bt) - f(\bx)| \leq C|\bt|^\delta.
\]

\begin{lemma} \label{lem:Fourier}Suppose that \( f: \mathbb{R}^d \to \mathbb{C} \) is an almost periodic function such that $\omega_1(f, \delta) \le C\delta^\eta$ for some constants $C>0$ and $\eta>0$ and for all $\delta >0$.
Then
\begin{equation}\label{eq:Holder}
\forall \bbeta \in \mathbb{R}^d\setminus \{\bm{0}\}, \quad |\widehat{f}(\bbeta)| \le \frac{C}{2}  \frac{1}{|\bbeta|^\eta}.
\end{equation}
\end{lemma}

\begin{proof}
Starting from the equalities  $$
\widehat{f}(\bbeta) = \mathcal{M}(f(\bx) e^{-2\pi i\bbeta^T  \bx})= \mathcal{M}(f(\bx+\btau) e^{-2\pi i\bbeta^T (\bx+\btau)})
$$
which hold for all $\bbeta \in \mathbb{R}^d$ and all $\btau\in \mathbb{R}^d$, 
we get
\[
\widehat{f}(\bbeta) \big(e^{2\pi i \bbeta^T \btau} -1\big)= \mathcal{M}\Big(\big[f(\bx+\btau) - f(\bx)\big] e^{-2\pi i\bbeta^T \bx}\Big).
\]
It follows that
\[
|\widehat{f}(\bbeta)| \cdot |e^{2\pi i \bbeta^T \btau}-1 | \leq \mathcal{M}(|f(\bx+\btau) - f(\bx)|).
\]
Notice that \( \btau \) in \( \mathbb{R}^d \) is arbitrary. For any
\( \bbeta \in \mathbb{R}^d\setminus \{\bm{0}\} \) with, take
\(
\btau = \pi \bbeta/|\bbeta|^2  
\)
(the vector in the direction of \( \bbeta \) with length \( (2|\bbeta|)^{-1} \)), so that
\(
\bbeta^T \btau = 1/2.
\)
Thus,
\(
e^{2\pi \bbeta^T \btau} = -1
\)
and we can conclude.
\end{proof}

The estimate \eqref{eq:Holder} is not effective for small $\bbeta$. However there are almost periodic functions having $0$ as accumulation point of its spectrum. For such functions,  the estimate \eqref{eq:Holder} does not supply any useful information about the  Fourier coefficients $\widehat{f}(\bbeta)$ for $\bbeta$ close to $\bm{0}$. In order to obtain some regularity of the parent function $F$  of $f$ basing on  the estimate \eqref{eq:Holder}, we need some extra condition on \( f \), some discreteness of its Bohr spectrum, 
as we have introduced by \eqref{eq:D-cond}.

Another argument is of Sobolev type.

\subsection{Density of \texorpdfstring{$\mathbf{P}(\bbZ^n)$}{} in \texorpdfstring{$\bbR^d$}{}}

As we have seen above, the spectrum $\sigma_f$ of any quasi-periodic function in ${\rm QP}_{\bP}(\mathbb{R}^d)$ is contained in 
$\bP(\mathbb{Z}^n)$. It is possible that $\bP(\mathbb{Z}^n)$ is dense in $\mathbb{R}^d$. 
The following is the condition for 
$\bP(\mathbb{Z}^n)$ to be dense in $\mathbb{R}^d$ when $d=1$.

\begin{lemma}\label{lem:density} The group $p_1\mathbb{Z} + \dots + {p}_n\mathbb{Z}$ is dense in $\mathbb{R}$ if and only if $\tfrac{p_i}{p_j}$ is irrational for some couple $(i, j)$. 
\end{lemma}
\begin{proof}
Assume that all $\tfrac{p_i}{p_j}$ are rational, equivalently $r_i:= \frac{p_i}{p_1}$ ($1\le i\le n$) are rational. Assume $r_i = \frac{s_i}{t_i}$
with integers $s_i, t_i$. We can write 
$r_i =\frac{s_i'}{t}$ with an integer $s_i'$ and $t=t_1t_2\cdots t_n$. It follows that $p_i= \frac{p_1}{t} s_i'$ and then
$$
     p_1\mathbb{Z} + \dots + {p}_n\mathbb{Z} 
     \subset  \frac{p_1}{t} \mathbb{Z}
$$
which is not dense in $\mathbb{R}$.

Now assume that $\alpha :=\frac{p_i}{p_j}$ is irrational for some couple $(i, j)$. Then $p_i\mathbb{Z}+p_j\mathbb{Z}$ is dense in $\mathbb{R}$ because 
$$
   p_i\mathbb{Z}+p_j\mathbb{Z}
   =p_j (\alpha \mathbb{Z}+\mathbb{Z})
$$
and $\alpha \mathbb{Z}+\mathbb{Z}$ is dense in $\mathbb{R}$, which is well known. It follows immediately the density of
$p_1\mathbb{Z} + \dots + {p}_n\mathbb{Z}$.
\end{proof}

\subsection{Smooth quasi-periodic function without smooth parent function}

Let $\phi=\frac{1+\sqrt{5}}{2}$ be the Golden ratio, which is the largest root of $x^2 -x-1=0$. Let $\phi'=\frac{1-\sqrt{5}}{2}$ the algebraic conjugate of $\phi$.  Let 
$$
\bm{A}=\begin{pmatrix}
    1 & \phi\\ 1 &\phi'
\end{pmatrix}.
$$
Consider the lattice
$$
\mathcal{L}: \bm{A}(\mathbb{Z}^2)=\{(m+n\phi, m+\phi' n):  m, n\in \mathbb{Z}\}.
$$
Take a bounded interval $W$ with non-empty interior. Define
\begin{eqnarray*}
    \Lambda:&=& \{m+\phi n: \ m,n\in \mathbb{Z}; m+n\phi'\in W \}\\
    \Lambda':&=& \{m+\phi' n: m,n\in \mathbb{Z}; m+n\phi'\in W \}\\
    \mathcal{B}:&=& \{ (m,n)\in \mathbb{Z}^2: m+n\phi'\in W \}.
\end{eqnarray*}
The set $\Lambda (\subset \mathbb{R})$ is a model set of Meyer \cite{Meyer1972}. In particular, it is uniformly discrete  and   relatively dense. It follows that for large $L>0$ and for all $a\in \mathbb{R}$ we have
\begin{equation}\label{eq:Model}
C L \le \#(\Lambda\cap [a, a+L])\le C' L
\end{equation}
for some constants $C>0$ and $C'>0$. On the other hand, the set $\Lambda'$ is contained in the bounded interval $W$ and it is dense in $W$.
Notice that $\mathcal{B}$ is a band in the lattice $\mathcal{L}$. The set $\Lambda$ and $\Lambda'$ are nothing but the projections of $\mathcal{B}$ to the first coordinate (the physical space) and the second coordinate (the internal space). 

\begin{lemma} \label{lem:golden} There exists a constant $C>1$ such that 
$$
\forall (m, n)\in \mathcal{B}, \quad \frac{1}{C} \sqrt{m^2+n^2}\le |m+n\phi|\le C \sqrt{m^2+n^2}.
$$
\end{lemma}
\begin{proof} Since the matrix $A$ is invertible, we have
$$
\frac{1}{C'} \sqrt{m^2+n^2}
\le \sqrt{|m+n\phi|^2 +|m+\phi'n|^2} \le C' \sqrt{m^2+n^2}
$$
 for all $(m, n)\in \mathbb{Z}^2$ and for some constant $C'>1$. We can conclude by using the facts that $|m+\phi'n|^2$ is bounded on $\mathcal{B}$ and $\inf_{(m,n) \in \mathcal{B}\setminus \{(0,0)\}}|m+\phi n|>0$ (recall that $\Lambda$ is uniformly discrete). 
    \end{proof}

Let us consider the quasi-periodic function
\begin{equation}\label{eq:example}
f(x) = \sum_{(n, m)\in \mathcal{B}\setminus \{(0,0)\} } \frac{1}{|n+\phi m|^{3/2}} e^{2\pi i (n+\phi' m) x}.
\end{equation}

\begin{proposition}
The above function $f$ is a quasi-periodic function of class $C^\infty$
and its derivatives $f^{(k)}$, $k\ge 1$, are also quasi-periodic. But its parent function is not in $C^1(\mathbb{T}^2)$. 
\end{proposition}

\begin{proof}
The key to prove the $C^\infty$ regularity is the convergence of the following series
$$
\sum_{(n, m)\in \mathcal{B}\setminus \{(0,0)\} } \frac{1}{|n+\phi m|^{3/2}} <\infty. 
$$
Indeed, the above series is a sum over $\Lambda$. We order $\Lambda \cap [0,+\infty)$ as an increasing sequence $(\lambda_k)$. As $\Lambda$ is uniformly discrete, $(\lambda_k)$ increases at least linearly, namely $\lambda_k \ge c k$ for some $c>0$. By the way, we remark that it increases also at most linearly, i.e. $\lambda_k\le b k$ for some $b>0$ because of the relative density of $\Lambda$. But, for the moment, we don't need the fact $\lambda_k\le b k$. So, up to a multiplicative constant, the sum over $\Lambda\cap[0, +\infty)$ is bounded by $\sum_{k\ge 1} \frac{1}{k^{3/2}}<\infty$. In the same way, we can prove that the sum over$\Lambda \cap(-\infty, 0)$ is finite too.

Now recall that $|n+\phi'm|$ is bounded for all $(n, m)\in \mathcal{L}^*$. Therefore, we can claim immediately that  $f$ is $k$-th differentiable for all $k\ge 1$ and we have
$$
f^{(k)}(x) = \sum_{(n, m)\in \Lambda_*\setminus \{(0,0)\} } \frac{(2\pi i (n+\phi'm)^k}{|n+\phi m|^{3/2}} e^{2\pi i (n+\phi' m) x}.
$$

The parent function $F$ of $f$ is developed into its Fourier series as follows
$$
F(y_1,y_2)  
= \sum_{(n, m)\in \mathcal{B}\setminus \{(0,0)\} } \frac{1}{|n+\phi m|^{3/2}} e^{2\pi i (n y_1+m y_2)}.
$$
If the partial derivatives $\frac{\partial F}{\partial y_1}$ and $\frac{\partial F}{\partial y_2}$ exist and are continuous (or in $L^2(\mathbb{T}^2)$), by Parseval identity we would have
$$
\sum_{(n, m)\in \mathcal{B}\setminus \{(0,0)\} } \frac{m^2+ n^2}{|n+\phi m|^3} <\infty. 
$$
But this series diverges, a contradiction. Indeed, by Lemma \ref{lem:golden}, the above series is comparable with  $\sum_{(n, m)\in \mathcal{B}\setminus \{(0,0)\} } \frac{1}{|n+\phi m|}=\infty$ (using the above remark: $\lambda_k \le b k$).
\end{proof}

We can use the above method to construct smooth quasi-periodic functions on $\mathbb{R}^d$ having non smooth parent functions, by using the cut-projection scheme of Meyer.  

The spectrum of the function $f$ constructed above is dense in the interval $W$. 
Now consider the following quasi-periodic function
\begin{equation}\label{eq:example2}
g_r(x) = \sum_{(n, m)\in \mathcal{B}\setminus \{(0,0)\} } \frac{1}{|n+\phi m|^{r}} e^{2\pi i (n+\phi m) x}.
\end{equation}
The spectrum of $g_r$ is the Meyer set $\Lambda\setminus\{(0,0)\}$.
This function $g_r$ has some degree of regularity, depending on  $r$.  Its parent function also has some degree of regularity.

\medskip

\subsection{Sobolev-Besicovitch spaces and regularity of parent functions} 

The Hausdorff-Young inequality shows that the Besicovitch norm $\|f\|_q$ 
of a Besicovitch function $f\in B^q_{ap}$ is comparable to some extent with the $\ell^{q'}$-norm of the Fourier coefficients of $f$. 
So, it is natural to introduce the following Sobolev-Besicovitch spaces through Fourier transforms (cf.~\cite{Besicovitch1932}).

For $1\le q<\infty$ and $s>0$, we define the Sobolev-Besicovitch space $H^{s,q}_{ap}(\mathbb{R}^d)$ by
$$
H^{s,q}_{ap}(\mathbb{R}^d)
=\left\{f \in B^{1}_{ap}(\mathbb{R}^d) : \|f\|_{H^s_{ap}}^{q'}=\sum_{\blam\in \sigma_f} (1+|\blam|^{2})|^{\frac{sq'}{2}}|\widehat{f}(\blam)|^{q'})<\infty 
\right\}
$$
where $q'=\frac{q}{q-1}$ is the conjugate of $q$.

\begin{theorem}\label{thm:SB-embedding}
    Let $1<q<\infty$. Let $f \in QP_{\bP}(\mathbb{R}^d)$ be a quasi-periodic function
    with parent function $F \in C(\mathbb{T}^n)$. Suppose  
    \begin{itemize}
        \item [(\,i\,)]
     $f\in H^{s, q}_{ap}(\mathbb{R}^d)$ for some $s>0$ such that $s-\frac{n}{q}>1$. 
    \item[(ii)]  $\bP \bk\|\ge C \|\bk\|$ for all $\bP \bk \in \sigma_f$ where $C$ is a constant.
    \end{itemize}
    Then $F\in C^m(\mathbb{T}^n)$ for all integer $m<s-\frac{n}{q}$.
\end{theorem}

\begin{proof}
It is well known that for any multi-index $\balpha =(\alpha_1, \cdots, \alpha_n)$ with $\sum_{j=1}^n \alpha_j \le m$, 
we have $|\bk^{\balpha}|=O(|\bk|^m)$. Therefore
we have only to show that $\sum |\bk|^m  |\widehat{F}(\bk)|<\infty$. Indeed, 
let $L_*=\{\bk \in \mathbb{Z}^n: \bP \bk \in \sigma_f\}$ which is the support of the Fourier transform of the parent function $F$. From
$$
\sum_{\bk\in \mathbb{Z}^n} |\bk|^m |\widehat{F}(\bk)|
= \sum_{\bk\in L_*} |\bk|^{m-s }|\cdot |\bk|^{s}\widehat{f}(\bP \bk)|,
$$
using the hypothesis $|\bP \bk|\ge C|\bk|$
for $\bP \bk\in \sigma_f$ and the H\"{o}lder inequality, we get 
$$
\sum_{\bk\in \mathbb{Z}^n} |\bk|^m |\widehat{F}(\bk)|
\le \frac{1}{C^s}\left(\sum_{\bk\in L_*} \frac{1}{|\bk|^{(s-m)q}} \right)^{1/q}
\left(\sum_{\blam\in \sigma_f} \Big(|\blam|^{sq'}|\widehat{f}(\blam)|\Big)^{q'}\right)^{1/q'}.
$$
The last factor is bounded by $\|f\|_{H^{s,q}_{ap}}$, which is finite. The sum of $\frac{1}{|\bk|^{(s-m)q}}$ is also finite, because $(s-m)q>n$.   
\end{proof}

When $q=1$ (i.e. $q'=\infty$), the conclusion still holds if the condition (i) on $f$ is replaced by $\widehat{f}(\blam) = O(|\blam|^{-s})$. Such a condition is really satisfied, see Lemma \ref{lem:Fourier}. See also Theorem \ref{thm:regularity}  

\subsection{Proof of Theorem \ref{thm:regularity}} 
For $\balpha:=(\alpha_1, \cdots, \alpha_n)^T \in \mathbb{Z}^n$, 
let $\beta = \bP\balpha \in \mathbb{R}$ where $\bP=(p_1, \cdots,p_n)$ is a $1\times n$ matrix and  where $\{p_1, \cdots, p_n\}$ is an integral basis of the $\mathbb{Z}$-module generated by $\sigma_f$. Then the parent function $F\in C(\mathbb{T}^n)$ is determined by $F(p_1 x, \cdots, p_n x)=f(x)$, namely $F(\bP^Tx)=f(x)$.

Since \( f\in C^{r, \eta}(\mathbb{R})\) and $\widehat{f^{(r)}} (\beta) = (2\pi i \beta)^r \widehat{f}(\beta)$, we have
\[
\widehat{f}(\beta) = \frac{1}{(2\pi i \beta)^r} \widehat{f^{(r)}}(\beta)= O\left(\frac{1}{|\beta|^{r+\eta}}\right)
\]
where  we have applied Lemma \ref{lem:Fourier} for the last equality. 
Since $ \widehat{F}(\balpha)=\widehat{f}(\bP\balpha)$, it follows that 
\[
 \forall \balpha \in \Sigma_F \ {\rm  with}\ \balpha\not=\bm{0}, \ \quad \widehat{F}(\balpha) = O\left(\frac{1}{|\bP \balpha|^{r+\eta}}\right) = O\left(\frac{1}{|\balpha|^{r+\eta}}\right),
\]
where for the last equality the assumption \eqref{eq:D-cond} is used. Recall that $|\balpha|$ denotes the Euclidean norm of $\balpha$.

Let $\bs=(s_1, \cdots, s_n)$ be a multiple index. 
Formally, we have the Bohr-Fourier series
\[
\frac{\partial^{\bs} F}{\partial x_1^{s_1}\cdots \partial x_n^{s_n}}(x) \ \  \sim \ \ \ (2\pi i)^{s_1+\cdots + s_n}
\sum_{\balpha\in  \sigma_F} \alpha_1^{s_1}\cdots \alpha_n^{s_n} \widehat{F}(\alpha) e^{2\pi i (\alpha_1 x_1 +\cdots + \alpha_n x_n)}.
\]
To conclude, it suffices to ensure that the last series converges absolutely. Indeed, 
\[\sum_{\balpha\in  \sigma_F} |\alpha_1^{s_1}\cdots \alpha_n^{s_n}| |\widehat{F}(\alpha)| \le 
  \sum_{\balpha\in  \sigma_F} \frac{1}{|\balpha|^{r+\eta -(s_1+\cdots s_n)}} <\infty
\]
if $r+\eta -(s_1+\cdots +s_n) >n$. As $\eta >0$, the finiteness of the last series is indeed ensured for all $(s_1, \cdots,s_n)$ such that $s_1+\cdots+s_n \le r-n$.

\section*{Acknowledgments}
	A. Fan was supported by the NSFC (No.~12231013), K. Jiang was supported by the National Key R\&D Program of China (No.~2023YFA1008802), P. Zhang was supported by the NSFC (No.~12288101).

\bibliographystyle{siamplain}
\bibliography{references}
\end{document}